\documentclass[%
 print,
superscriptaddress,
 amsmath,amssymb,
aps,
]{revtex4-2}

\usepackage{graphicx}
\usepackage{dcolumn}
\usepackage{bm}

\usepackage{algorithm}
\usepackage{algpseudocode}
\usepackage{multirow}

\renewcommand{\b}[1]{\boldsymbol{#1}}
\newcommand{\T}[0]{^\mathrm{T}}

\begin{document}

\title{SINDy on slow manifolds}

\author{Diemen Delgado-Cano}
\affiliation{Department of Mechanical Engineering, Universidad de Chile, Beauchef 851, Santiago, Chile}

\author{Erick Kracht}%
\affiliation{Department of Mechanical Engineering, Universidad de Chile, Beauchef 851, Santiago, Chile}

\author{Urban Fasel}%
\affiliation{Department of Aeronautics, Imperial College London, SW7 2AZ, United Kingdom}

\author{Benjamin Herrmann}%
 \email{benjamin.herrmann@uc.cl}
\affiliation{Department of Mechanical and Metallurgical Engineering, Pontificia Universidad Católica de Chile, Av. Vicuña Mackenna 4860, Santiago, Chile}
\affiliation{Department of Hydraulic and Environmental Engineering, Pontificia Universidad Católica de Chile, Av. Vicuña Mackenna 4860, Santiago, Chile}


\begin{abstract}
The sparse identification of nonlinear dynamics (SINDy) has been established as an effective method to learn interpretable models of dynamical systems from data. However, for high-dimensional slow-fast dynamical systems, the regression problem becomes simultaneously computationally intractable and ill-conditioned. Although, in principle, modeling only the dynamics evolving on the underlying slow manifold addresses both of these challenges, the truncated fast variables have to be compensated by including higher-order nonlinearities as candidate terms for the model, leading to an explosive growth in the size of the SINDy library. In this work, we develop a SINDy variant that is able to robustly and efficiently identify slow-fast dynamics in two steps: (i) identify the slow manifold, that is, an algebraic equation for the fast variables as functions of the slow ones, and (ii) learn a model for the dynamics of the slow variables restricted to the manifold. Critically, the equation learned in (i) is leveraged to build a \textit{manifold-informed} function library for (ii) that contains only essential higher-order nonlinearites as candidate terms. Rather than containing all monomials of up to a certain degree, the resulting custom library is a sparse subset of the latter that is tailored to the specific problem at hand. The approach is demonstrated on numerical examples of a snap-through buckling beam and the flow over a NACA 0012 airfoil. We find that our method significantly reduces both the condition number and the size of the SINDy library, thus enabling accurate identification of the dynamics on slow manifolds.
\end{abstract}

\maketitle


\section{\label{sec:intro}Introduction}

Slow-fast dynamical systems are pervasive in science and engineering, arising in contexts such as fluid mechanics, structural dynamics, and climate modeling. These systems are characterized by high-dimensional state spaces and with trajectories that, due to timescale separation, often collapse onto lower-dimensional slow manifolds, where the effective dynamics evolve. These manifolds include stable, unstable, and center manifolds from classic nonlinear dynamics theory~\citep{GyHbook}, as well as the more recently established spectral submanifolds (SSMs), defined as the smoothest nonlinear continuation of a nonresonant spectral subspace of a hyperbolic steady state~\citep{haller2016nd}. The dynamics restricted to each of these manifolds provide a reduced-order representation of the long-term behavior of the system. Therefore, understanding and modeling the effective dynamics on slow manifolds is crucial for efficient simulation, control, and prediction of high-dimensional nonlinear systems.

Recent advances in data-driven modeling have enabled the discovery of low-order representations of high-dimensional nonlinear dynamical systems directly from measurements. The dynamic mode decomposition (DMD)~\citep{schmid2010jfm,dmdbook}, along with its many variants~\citep{herrmann2021jfm,baddoo2022prsa,baddoo2023prsa,schmid2022arfm}, are effective analysis tools to extract spatio-temporal patterns from data. However, they are not suitable to build low-order models that capture inherently nonlinear dynamics. Operator inference (OpInf)~\citep{peherstorfer2016cmame,qian2020physd,kramer2024arfm} is a prominent approach that learns low-order nonlinear models by projecting state data onto a basis of proper orthogonal decomposition (POD) modes~\citep{berkooz1993arfm} and then identifies the dynamics of the modal coefficients using polynomial regression. Typically, when using OpInf, the maximum degree of the monomials used for regression is selected based on physical knowledge of the structure of the underlying governing equations, resulting in the identified model being an approximation of a POD-Galerkin reduced order model (ROM) of the system~\citep{peherstorfer2016cmame}. Recently, a related approach was developed to simultaneously identify oblique projection operators and a Petrov-Galerkin ROM from data~\citep{padovan2024siads}. Also recently, data-driven SSM identification has emerged as a powerful tool to extract nonlinear dynamics evolving on low-dimensional SSMs from data~\citep{cenedese2022natcomms,axaas2023nd_a}. The approach is based on the rigorous mathematical foundations that underpin SSM theory to allow for the identification of the effective degrees of freedom and governing equations on the SSM, offering a path toward interpretable and efficient low-order models for nonlinear systems, although limited to manifolds that are anchored to a hyperbolic steady state. Nonetheless, the identified models have been shown to capture chaotic dynamics~\citep{liu2024chaos} and the forced response of mechanical systems~\citep{cenedese2022ptrsa}.

The sparse identification of nonlinear dynamics (SINDy)~\citep{brunton2016pnas}, has been established as a powerful framework to learn interpretable models of dynamical systems from data. The method leverages sparsity promoting regression to fit a sparse set of terms from a predefined function library to the observed dynamics, thus balancing model complexity and accuracy to avoid overfitting~\citep{brunton2016pnas}. Since its introduction, several innovations have been developed~\citep{kaiser2018prsoc,fasel2022prsa,lemus2024nd,brunton2025arcras} that extend the capabilities of the original SINDy algorithm, many have been implemented in the open source software package PySINDy~\citep{desilva2020joss,kaptanoglu2022joss}, and some have been recently benchmarked~\citep{kaptanoglu2023nd} using the dysts standardized database of chaotic systems introduced by~\cite{gilpin2021neurips}. Although the original SINDy formulation was limited to systems with only a few state variables, its combination with dimensionality reduction techniques, such as the use of linear embeddings given by POD or DMD modes~\citep{loiseau2018jfm_a,loiseau2020dtcfd}, or deep autoencoders~\citep{champion2019pnas,fukami2021jfm,fries2022cmame,callaham2022jfm,bakarji2023prsa,conti2023cmame}, has become common practice to enable modeling of high-dimensional dynamics evolving on low-dimensional latent spaces. It should be noted that SINDy has been widely adopted by a large number of practitioners, mainly due to its simple philosophy, ease of implementation, and effectiveness.

For the specific case of slow-fast systems, as we detail in section~\S\ref{sec:background2}, fast variables may be expressed in terms of the slow ones via a manifold equation and, therefore, the dynamics are fully driven by the slow variables. The implication for both SINDy and OpInf is that the effect of the fast variables needs to be accounted for by including higher-order nonlinearities in the model ansatz for the slow dynamics. If this is not taken into account, the underlying regression problems become ill-conditioned, leading to inaccurate models with poor predictive performance. We discuss the root cause of this ill-conditioning in section~\S\ref{sec:methods2}. The work in~\citep{callaham2022jfm} addressed this in the context of SINDy for systems with nonlinear correlations between state variables by separately identifying a manifold equation and the dynamics of the driving variables restricted to said manifold. Similarly, recent extensions of the OpInf framework are able to learn low-order nonlinear dynamics evolving on quadratic~\citep{geelen2023cmame} and more general polynomial manifolds~\citep{geelen2024chaos}. Moreover, data-driven SSM identification naturally addresses the learning of both manifold and dynamics by construction~\citep{cenedese2022natcomms}. The more recent version of the method separates these two learning tasks and performs them sequentially, drastically reducing the computational complexity of the underlying optimization problem to improve the scaling properties with the dimension of the system~\citep{axaas2023nd_a}. Even though SINDy, OpInf, and SSM identification are all able to learn dynamics restricted to a slow manifold, they all rely on polynomial regression onto a dense function library that may need to include monomials of very high degree to properly account for strong nonlinearities due to the effect of the fast variables.

In this work, we propose a novel extension of SINDy tailored for slow-fast systems. Similarly to previous methods, our approach consists of two steps: (i) identifying an algebraic equation that defines the slow manifold and (ii) learning a model for the slow dynamics restricted to the manifold. Importantly, we construct a manifold-informed function library that incorporates only the necessary higher-order nonlinearities, significantly reducing the size of the SINDy library and improving numerical conditioning. This targeted function library circumvents the explosion in the amount of candidate terms typically encountered in slow-fast systems, allowing for more accurate and efficient identification of reduced-order models.

The remainder of the paper is organized as follows. Theoretical background on SINDy and slow-fast systems is covered in~\S\ref{sec:background}. Our proposed SINDy extension for dynamics on slow manifolds is formulated in~\S\ref{sec:methods}. The method is demonstrated on numerical examples and its performance is assessed and discussed in~\S\ref{sec:results}. Our conclusions are offered in~\S\ref{sec:conclusions}.

\section{\label{sec:background}Background}

In this section we provide an introduction to SINDy in the context of high-dimensional systems and a brief review of slow-fast dynamics.

\subsection{\label{sec:background1}SINDy for high-dimensional systems}

In this exposition, we take as a starting point a high-dimensional nonlinear dynamical system of the form
\begin{equation}
    \b{\dot{q}}=\b{\tilde{F}}(\b{q}),\label{sys_full}
\end{equation}
where the overdot denotes time-differentiation, $\b{q}\in\mathbb{R}^{n_q}$ is the state of the system, and $\b\tilde{{F}}: \mathbb{R}^{n_q} \rightarrow \mathbb{R}^{n_q}$ its dynamics that are assumed to be unknown. Given a set of $m$  measurements of the state $\lbrace\b{q}_j\rbrace$ for $j=1,\dots,m$, typically acquired from a family of trajectories of interest, we wish to find an interpretable and predictive model for the high-dimensional dynamics. This can be achieved combining dimensionality reduction techniques and SINDy, as described below.

It is often the case that, over the state space region of interest, the dynamics evolve on a much lower dimensional latent space $\b{z}\in\mathbb{R}^{n_z}$, with $n_z \ll n_q$, that we can map to and from with the transformations
\begin{equation}
    \b{z}=\mathcal{E}(\b{q}), \text{ and } \b{q}=\mathcal{D}(\b{z}),\label{autoencoder}
\end{equation}
referred to as encoder and decoder, respectively. Note that the dynamics of the latent representation are then related to the full state dynamics via $\b{\dot{z}}=\b{F}(\b{z})=\mathcal{E}(\b{\tilde{F}}(\mathcal{D}(\b{z})))$. The transformations in eq.~\eqref{autoencoder} can be learned directly from the measurements of $\b{q}$ using, for example, deep autoencoder neural networks~\cite{champion2019pnas,fukami2021jfm,callaham2022jfm,fries2022cmame,conti2023cmame}. Moreover, a linear embedding can be obtained from a modal analysis of the linearized dynamics~\citep{khoo2022jfm} or via linear dimensionality reduction techniques~\citep{loiseau2018jfm_a,loiseau2020dtcfd}. In the linear embedding scenario, the encoding and decoding transformations in eq.~\eqref{autoencoder} become matrix multiplications 
\begin{equation}
    \b{z}=\b{W_z}\T\b{q}, \quad \text{and} \quad \b{q}=\b{V_z}\b{z},\label{modal}
\end{equation}
with $\b{W_z}$ and $\b{V_z}\in\mathbb{R}^{n_q \times n_z}$ being matrices containing bi-orthonormal sets of vectors such that $\b{W_z}\T \b{V_z}=\b{I}$. Hence, the latent dynamics become $\b{F}(\b{z})=\b{W_z}\T\b{\tilde{F}}(\b{V_z}(\b{z}))$, which is an (oblique) Petrov-Galerkin projection of the full dynamics~\citep{rowley2017arfm}. Note that this becomes an (orthogonal) Galerkin projection if $\b{W_z}=\b{V_z}$~\citep{rowley2017arfm}. Here we assume that, for all practical purposes, this linear or nonlinear dimensionality reduction step acts as a lossless compression, eliminating only redundant degrees of freedom in the original representation of the system. For the rest of this article, we continue the exposition in the setting of linear dimensionality reduction with multiplication with $\b{W_z}\T$ and $\b{V_z}$ replacing the action of the nonlinear transformations $\mathcal{E}(\cdot)$ and $\mathcal{D}(\cdot)$, respectively, although this is without loss of generality. 

The goal now shifts to identifying a model for the low-order dynamics $\b{F}$ to predict $\b{z}$ from which the full dimensional state $\b{q}$ can be reconstructed. SINDy assumes that the dynamics can be expressed as a linear combination of a few nonlinear functions selected among a much larger, and judiciously designed, set of candidate functions, as follows
\begin{equation}
    \b{\dot{z}}=\b{F}(\b{z})\approx \b{\Xi_z}\T\b{\theta_z}(\b{z}),\label{ansatz_z}
\end{equation}
where $\b{\theta_z} : \mathbb{R}^{n_z} \rightarrow \mathbb{R}^{\ell_z}$ contains the set of $\ell_z$ candidate scalar functions, and $\b{\Xi_z}\in\mathbb{R}^{\ell_z\times n_z}$ are the yet to be identified model coefficients.

Returning to our dataset built from state observations, we may assemble the data matrices
\begin{equation}
    \b{\dot{Z}} = \left[\b{\dot{z}}_1 \ \b{\dot{z}}_2 \ \cdots \ \b{\dot{z}}_m\right]\T  \in\mathbb{R}^{m\times n_z},
    \quad \text{and} \quad
    \b{\Theta_z}=\left[\b{\theta_z}(\b{z}_1) \ \b{\theta_z}(\b{z}_2) \ \cdots \ \b{\theta_z}(\b{z}_m)\right]\T \in\mathbb{R}^{m\times \ell_z},\label{data_z}
\end{equation}
where $\b{z}_j=\b{W_z}\T \b{q}_j$ and the time derivatives $\b{\dot{z}}_j=\b{\dot{z}}(t_j)$ may be approximated from sequential data for the latent variable, for example, via finite differences. This approximation might be problematic if dealing with noisy data, and more recent variants use a weak form formulation to avoid differentiation~\citep{schaeffer2017pre,messenger2021mms}. We may now formulate an optimization problem to identify $\b{\Xi_z}$, as follows
\begin{equation}
  \b{\Xi_z} = \underset{\b{\Xi_z'}}{\mathrm{argmin}} \ \| \b{\dot{Z}}  -\b{\Theta_z} {\b{\Xi_z'}} \|_F^2 + \lambda_z \|\b{\Xi_z'}\|_0,\label{sindy_opt}
\end{equation}
where the first term in the cost function accounts for how well the model fits the data, and the second term, weighted by the hyperparameter $\lambda$, penalizes the amount of non-zero entries in the coefficient matrix, that is, the amount of terms in the identified model. Finding an exact solution to~\eqref{sindy_opt} requires performing a brute force search over all possible model structures, which is intractable when the number of candidate functions is large. In practice, sparse regression techniques bypass the combinatorial search by only approximating the optimal solution to~\eqref{sindy_opt} by, either solving a convex relaxation of the problem, such as the LASSO~\citep{tibshirani1996jrsssb}, or using greedy algorithms, such as the sequentially thresholded least squares (STLS) introduced in the original SINDy article~\citep{brunton2016pnas}. In this work we use the latter due to its simplicity and proven effectiveness in practice.

\subsection{\label{sec:background2}Dynamics restricted to a slow manifold}

In this work, we loosely define a stable (unstable) slow manifold as an invariant manifold over which trajectories evolve slowly and towards which trajectories are attracted exponentially fast in forward (backward) time, at least locally in a region of interest in state space. This includes stable, center, and unstable manifolds, as well as SSMs, of a system at a stationary state, such as a fixed point or a periodic or quasiperiodic orbit. For an interesting discussion on more rigorous definitions of slow manifolds, see~\citep{lorenz1992jas}.

Slow manifolds emerge naturally in slow-fast systems where there is a separation of time scales. In that scenario, we may express the dynamics for its slow and fast components as follows
\begin{equation}
    \b{\dot{x}}=\b{f}(\b{x},\b{y}), \quad \text{and} \quad
    \b{\dot{y}}=\b{g}(\b{x},\b{y}), \label{sys_slowfast}
\end{equation}
where $\b{x}\in\mathbb{R}^{n_x}$ are $\b{y}\in\mathbb{R}^{n_y}$ are the slow and fast variables, respectively, and $\b{f}$ and $\b{g}$ their dynamics. For trajectories on the slow manifold, the fast variables are \emph{slaved} to the slow ones such that, at least locally, they may be expressed as a graph, yielding an algebraic equation of the form
\begin{equation}
    \b{y}=\b{h}(\b{x}).\label{manifold}
\end{equation}
As a consequence, the dynamics restricted to the slow manifold become
\begin{equation}
    \b{\dot{x}}=\b{f}(\b{x},\b{h}(\b{x})),\label{sys_on_manifold}
\end{equation}
which is a closed-form differential equation describing the evolution of $\b{x}$. Therefore, the dynamics of the full state are described by the system of differential-algebraic equations (DAEs) formed by equations~\eqref{manifold} and~\eqref{sys_on_manifold}. Moreover, the dynamics on an invariant manifold are tangent to the underlying vector field, requiring that
\begin{equation}
    \b{g}(\b{x},\b{h}(\b{x}))=\b{Dh}(\b{x}) \b{f}(\b{x},\b{h}(\b{x})),\label{tangency}
\end{equation}
where $\b{Dh}(\b{x})$ is the Jacobian of $\b{h}$ evaluated at $\b{x}$. This equation, known as the tangency condition, is obtained by differentiating eq.~\eqref{manifold}, applying the chain rule, and substituting $\b{\dot{x}}$ and $\b{\dot{y}}$ evaluated on the manifold from eqs.~\eqref{sys_slowfast}. When the governing equations are available, eq.~\eqref{tangency} is typically leveraged to build analytical approximations to the invariant manifold, for example via Taylor series expansions~\cite{GyHbook}, as shown in Appendix~\ref{app1}.

For completeness, we remark that, starting from the familiar form of system~\eqref{sys_full}, equations~\eqref{sys_slowfast} may be obtained by first finding linear embeddings such that
\begin{equation}
    \b{x}=\b{W_x}\T\b{q}, \quad  \b{y}=\b{W_y}\T\b{q}, \quad \text{and} \quad \b{q}=\b{V_x}\b{x} + \b{V_y}\b{y},\label{modal_slowfast}
\end{equation}
where $\b{W_x}\in\mathbb{R}^{n_q \times n_x}$ and $\b{W_y}\in\mathbb{R}^{n_q \times n_y}$ provide the linear transformations onto the slow and fast subspaces defined by the columns of $\b{V_x}\in\mathbb{R}^{n_q \times n_x}$ and $\b{V_y}\in\mathbb{R}^{n_q \times n_y}$, respectively. Importantly, these matrices satisfy
\begin{equation}
   \b{W_x}\T\b{V_x}=\b{W_y}\T\b{V_y} = \b{I}, \quad  \text{and} \quad \b{W_y}\T\b{V_x}=\b{W_x}\T\b{V_y}=\b{0},\label{modal_slowfast_conditions}
\end{equation}
where we are allowing an abuse of notation, since the above identity and zero matrices arising from each equality may have different dimensions. These matrices may be obtained, for example, via modal analysis of an equilibrium point. Subsequently, we may arrive at eqs.~\eqref{sys_slowfast} by defining
\begin{equation*}
    \b{f}(\b{x},\b{y}) = \b{W_x}\T\b{\tilde{F}}(\b{V_x}\b{x} + \b{V_y}\b{y}) \quad \text{and} \quad
    \b{g}(\b{x},\b{y}) = \b{W_y}\T\b{\tilde{F}}(\b{V_x}\b{x} + \b{V_y}\b{y}).
\end{equation*}
An example of simple dynamical system exhibiting a slow manifold is presented in the next section.

\section{\label{sec:methods}SINDy on slow manifolds}

In this section we formalize the setup for the data-driven modeling problem addressed in this work. We explain the pitfalls of a straightforward application of SINDy in this setting, present our proposed method, and discuss the connection to other existing approaches.

\subsection{\label{sec:methods1}Problem setup}

We begin from the same starting point as in section~\S\ref{sec:background1}, with a set of observations of the state $\b{q}$ of an unknown high-dimensional nonlinear dynamical system of the form of eq.~\eqref{sys_full} that we wish to identify. Importantly, we include as an additional consideration the fact that the system to be identified is known to exhibit slow-fast dynamics. Moreover, we assume to have access to the matrices $\b{W_z}$ and $\b{V_z}$ that provide the linear encoding and decoding transformations in eq.~\eqref{modal} to map to and from the latent variable $\b{z}$ to eliminate redundant states. Furthermore, we assume to know a split of these mappings 
\begin{equation}
    \b{W_z}=\left[\b{W_x} \ \b{W_y}\right], \quad  \text{and} \quad
    \b{V_z}=\left[\b{V_x} \ \b{V_y}\right], \label{modal_slowfast_split}
\end{equation}
such that conditions~\eqref{modal_slowfast_conditions} are satisfied and the latent variable is split into its slow and fast components, $\b{z}=[\b{x}\T \ \b{y}\T]\T$, via eq.~\eqref{modal_slowfast}. Depending on the application, $\b{W_x}$, $\b{W_y}$, $\b{V_x}$ and $\b{V_y}$ may be computed using modal analysis of the original system, or via data-driven techniques~\cite{axaas2023nd_a}. Therefore, given these four matrices, we can map our data for the full state $\b{q}_j$ onto measurements of the slow and fast variables $\b{x}_j$ and $\b{y}_j$, respectively. From this data, the goal is now to identify an interpretable and predictive dynamical system model. Critically, since the dynamics of the underlying system evolve on a slow manifold, we should seek a model with a DAE structure. We present the consequences of not enforcing this structure in the following subsection. 

\subsection{\label{sec:methods2}The origin of ill-conditioning in manifold-agnostic identification}

Let us consider the direct application of the approach presented in section~\S\ref{sec:background1} to data coming from a slow-fast system. This approach completely disregards the existence of the underlying slow manifold, therefore we refer to it as a manifold-agnostic identification procedure. As previously explained, for a slow-fast system, the dynamics quickly collapse onto a slow manifold where the fast variables are slaved to the slow ones so that $\b{y}=\b{h}(\b{x})$. As a result, the library of candidate functions in the model ansatz~\eqref{ansatz_z} becomes $\b{\theta_z}(\b{z})=\b{\theta_z}(\b{x},\b{y})=\b{\theta_z}(\b{x},\b{h}(\b{x}))$, which most likely contains linearly dependent functions of $\b{x}$. This is certainly the case if, as is commonly done in practice, we use monomials for the candidate functions, since monomials of the entries in $\b{x}$ will be linearly dependent with those arising from the Taylor expansion approximation of $\b{h}(\b{x})$. In other words, when building a library that includes functions of $\b{y}$, these will be redundant with some, already included, functions of $\b{x}$. As a consequence, the data matrix $\b{\Theta_z}$, built as in eq.~\eqref{data_z}, will be rank deficient, leading to an ill-conditioned regression problem and, in turn, a system identification procedure that is sensitive to noise in the data and produces models with poor predictive capabilities.

\begin{figure}
    \centering
    \includegraphics[width=17.5cm]{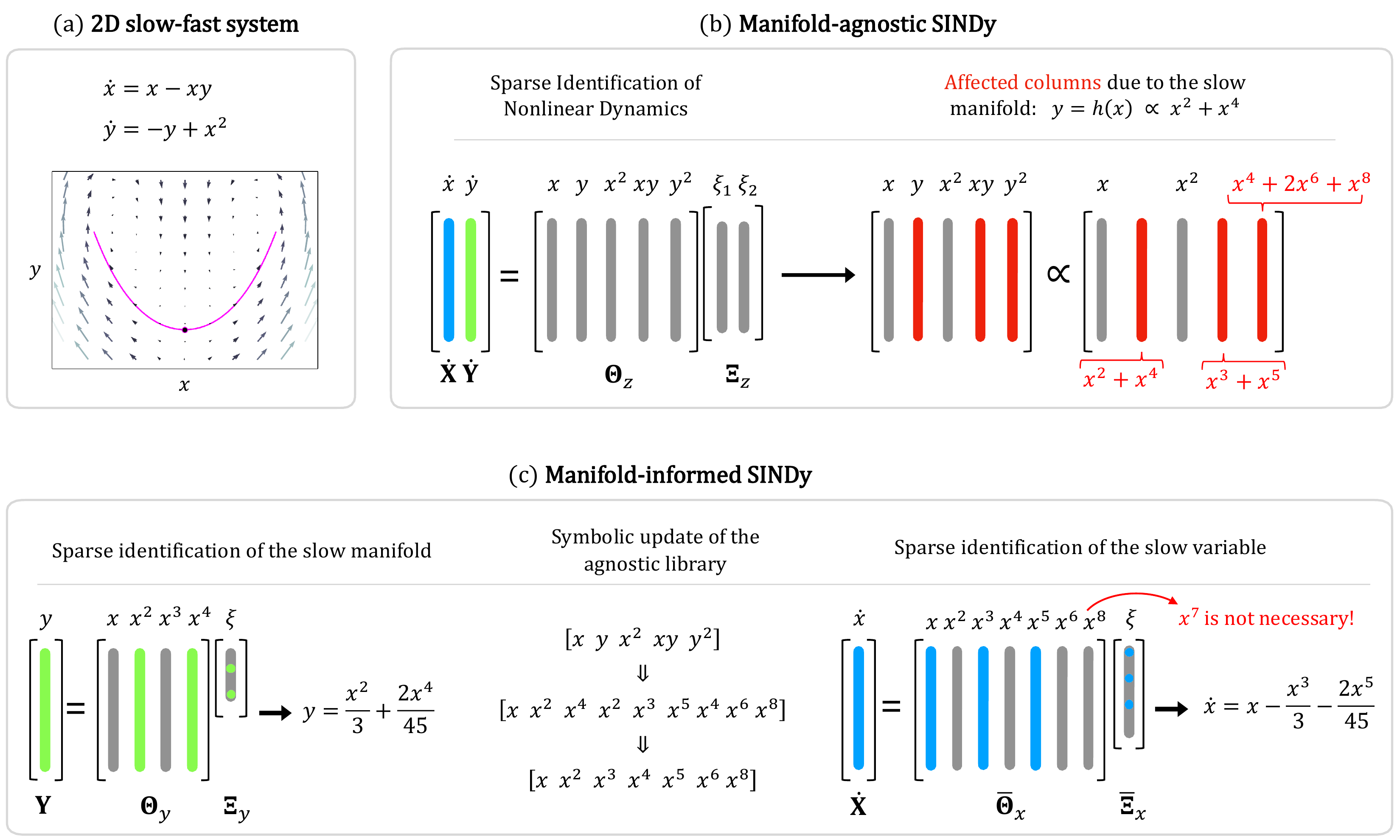}
    \caption{Comparison of the application of a manifold-agnostic library and a manifold-informed library. (a) A two-dimensional slow-fast system. (b) Application of SINDy with a manifold-agnostic library and the redundancy of its columns due to the presence of a slow manifold. (c) Identification of the slow manifold equation,  assembly of a manifold-informed library, and identification of the slow dynamics.}
    \label{fig:SINDyComparison}
\end{figure}

As a pedagogical example to illustrate this point, we consider the following low-dimensional dynamical system

\begin{subequations}
\begin{align}
    \dot{x}&=x-xy,\\
    \dot{y}&=-y+x^2.
\end{align}\label{eq:toy_system}
\end{subequations}
This system has an unstable manifold attached to its unstable equilibrium at $\b{z}=(x,y)=(0,0)$. Close to the origin, this manifold is well approximated by $y=h(x)=x^2/3 + 2x^4/45$, which can be obtained with a classic Taylor series approach, as detailed in Appendix~\ref{app1}. The phase plane for this system along with the manifold are shown in Fig.~\ref{fig:SINDyComparison}(\textit{a}). The straightforward application of SINDy considering observations of the full state on the manifold results in a rank-deficient function library data matrix. This happens because monomials of the state variables are used and, in view of the manifold equation, the term $x^2$ is linearly dependent with $y$, as shown in Fig.~\ref{fig:SINDyComparison}(\textit{b}).

The dynamics of the system restricted to the manifold are governed by the following DAE
\begin{subequations}
\begin{align}
    \dot{x}&=x-x^3/3-2x^5/45,\\
    y&=x^2/3+2x^4/45.\label{eq:toy_system DAE}
\end{align}
\end{subequations}
Recent approaches, described in the next subsection, avoid the conditioning difficulties discussed above by appropriately taking into account the DAE structure, as shown in Fig.~\ref{fig:SINDyComparison}(\textit{c}).

\subsection{Library size scaling in manifold-aware identification}

To enable accurate and robust identification of the dynamics of slow-fast systems, state-of-the-art methods, including SSM identification~\citep{cenedese2022natcomms,axaas2023nd_a} and recent variants of SINDy~\citep{callaham2022jfm} and OpInf~\citep{geelen2024chaos},  incorporate knowledge of the presence of a slow manifold into the respective regression problems. Because of this, we refer to this group of techniques as manifold-aware identification methods. Here we present the formulation for the manifold-aware version of SINDy~\citep{callaham2022jfm}, drawing the connections to and discussing the similarities with fast SSM identification~\citep{axaas2023nd_a} and OpInf~\citep{geelen2024chaos}. The underlying strategy is to divide the modeling effort into two sequential tasks: (i) identify an algebraic equation for the slow manifold, and (ii) identify the dynamics of the slow variables. Keep in mind that this is primarily enabled by the knowledge of the split between slow and fast dynamics. It is worth noting that, manifold-aware OpInf finds this split automatically, as it simultaneously identifies an embedding and the dynamics, although at an increased computational cost~\citep{geelen2024chaos}.

To begin with the first task, we assume the following ansatz for the slow manifold equation
\begin{equation}
    \b{y}=\b{h}(\b{x})\approx \b{\Xi_y}\T\b{\theta_y}(\b{x}),\label{ansatz_y}
\end{equation}
where $\b{\theta_y} : \mathbb{R}^{n_x} \rightarrow \mathbb{R}^{\ell_y}$ contains the set of $\ell_y$ candidate scalar functions, and $\b{\Xi_y}\in\mathbb{R}^{\ell_y\times n_y}$ are the yet to be identified model coefficients. Next, from our full state observations, we assemble the following data matrices
\begin{equation}
    \b{Y} = \left[\b{y}_1 \ \b{y}_2 \ \cdots \ \b{y}_m\right]\T  \in\mathbb{R}^{m\times n_y},
    \quad \text{and} \quad
    \b{\Theta_y}=\left[\b{\theta_y}(\b{x}_1) \ \b{\theta_y}(\b{x}_2) \ \cdots \ \b{\theta_y}(\b{x}_m)\right]\T \in\mathbb{R}^{m\times \ell_y},\label{data_y}
\end{equation}
recalling that our snapshots for the fast and slow variables are obtained from the full state as $\b{y}_j=\b{W_y}\T \b{q}_j$ and $\b{x}_j=\b{W_x}\T \b{q}_j$. Notice that there are no time derivatives in $\b{Y}$ since, in this instance, we are looking for an algebraic equation rather than an differential one. Now we identify $\b{\Xi_y}$ using the STLS algorithm to approximate the solution to
\begin{equation}
  \b{\Xi_y} = \underset{\b{\Xi_y'}}{\mathrm{argmin}} \ \| \b{Y}  -\b{\Theta_y} {\b{\Xi_y'}} \|_F^2 + \lambda_y \|\b{\Xi_y'}\|_0,\label{eq:cost_manifold}
\end{equation}
where $\lambda_y$ is the hyperparameter that promotes sparsity in the solution. The reasoning behind the use of a sparse regression for this step is that the function being learned is typically sparse in the space of candidate functions due to symmetries of the underlying manifold. On the other hand, at this step, SSM identification and OpInf consider similar regression problems, but without the sparsity penalty and including Tikhonov regularization for the latter, that are solved via least squares. In all these approaches, we typically consider as candidate functions all the monomials formed from the entries in $\b{x}$ that have up to a certain maximum degree $d_y$, the logic being to build a Taylor series expansion approximation of the underlying function $\b{h}(\b{x})$. The value specified for $d_y$ is usually selected via cross-validation. 

Once the slow manifold equation has been identified, we proceed with the identification of the dynamics for the slow variables governed by an equation of the form of eq.~\eqref{sys_on_manifold}. As before, we assume a linear combination of candidate nonlinear functions as an ansatz for the sought out function
\begin{equation}
    \b{\dot{x}}=\b{f}(\b{x},\b{h}(\b{x}))\approx \b{\Xi_x}\T\b{\theta_x}(\b{x}),\label{ansatz_x}
\end{equation}
where $\b{\theta_x} : \mathbb{R}^{n_x} \rightarrow \mathbb{R}^{\ell_x}$ contains the set of $\ell_x$ candidate scalar functions, and $\b{\Xi_x}\in\mathbb{R}^{\ell_x\times n_x}$ are model coefficients to be identified. Again, building towards a sparse regression problem, we assemble the following data matrices
\begin{equation}
    \b{\dot{X}} = \left[\b{\dot{x}}_1 \ \b{\dot{x}}_2 \ \cdots \ \b{\dot{x}}_m\right]\T  \in\mathbb{R}^{m\times n_x},
    \quad \text{and} \quad
    \b{\Theta_x}=\left[\b{\theta_x}(\b{x}_1) \ \b{\theta_x}(\b{x}_2) \ \cdots \ \b{\theta_x}(\b{x}_m)\right]\T \in\mathbb{R}^{m\times \ell_x},\label{data_x}
\end{equation}
where the time derivatives $\b{\dot{x}}_j=\b{\dot{x}}(t_j)$ are approximated from sequential data via finite differences. As may be expected, the identification now follows by using the STLS algorithm to approximate the solution to
\begin{equation}
  \b{\Xi_x} = \underset{\b{\Xi_x'}}{\mathrm{argmin}} \ \| \b{\dot{X}}  -\b{\Theta_x} {\b{\Xi_x'}} \|_F^2 + \lambda_x \|\b{\Xi_x'}\|_0,
\end{equation}
where $\lambda_x$ is the sparsity promoting hyperparameter. Again, for this step, SSM identification and OpInf solve a least squares problem to minimise the above cost function without the sparsity penalty term and with Tikhonov regularisation for the latter. 

Importantly, all the manifold-aware identification methods discussed usually assume a polynomial ansatz of a certain degree $d_x$. For a system with $n_x$ slow variables, the number of candidate terms in the dense polynomial library $\b{\theta_x}$ of degree $d_x$ can be computed as 
\begin{equation}
    \ell_x = \sum_{j=0}^{d_x}\frac{(n_x+j-1)!}{j!(n_x-1)!}.\label{eq:dense_library_size}
\end{equation}
Critically, eq.~\eqref{eq:dense_library_size} reveals that the library size grows drastically with the number of slow variables (dimension of the manifold) and the degree of nonlinearity of the slow dynamics. Moreover, if the underlying physical system is known to have a polynomial structure with degree $d_z$, then the maximum degree for the dynamics restricted to the slow manifold will be at most $d_x=d_z d_y$, where $d_y$ is the polynomial degree of the already identified manifold equation. Therefore, there is an explosion in library size with the degree of nonlinearity of the original system in physical space, and also with the curvature of the manifold.
In these scenarios, the library size becomes the computational bottleneck for manifold-aware identification methods.

\subsection{Proposed method: manifold-informed identification}

We propose a SINDy variant for slow-fast systems, schematically summarised in Fig.~\ref{fig:SINDyComparison}(\textit{c}), that combines knowledge of the physics of the system with the structure of the underlying manifold to build a leaner function library for the slow dynamics, containing significantly fewer terms than those in a dense polynomial library. The approach directly leverages the identified manifold equation to populate the SINDy library, and so we refer to the method as a manifold-informed method, and to the library as a manifold-informed library.

The new method follows the manifold-aware approach up to the identification of the manifold through the optimisation problem in eq.~\eqref{eq:cost_manifold}. A critical part of the proposed approach is the way in which we populate the manifold-informed ansatz for the slow dynamics
\begin{equation}
    \b{\dot{x}}=\b{f}(\b{x},\b{h}(\b{x}))\approx \b{\bar{\Xi}_x}\T\b{\bar{\theta}_x},(\b{x}),\label{ansatz_x}
\end{equation}
where $\b{\bar{\theta}_x} : \mathbb{R}^{n_x} \rightarrow \mathbb{R}^{\bar{\ell}_x}$ contains the set of $\bar{\ell}_x$ candidate scalar functions, and $\b{\bar{\Xi}_x}\in\mathbb{R}^{\bar{\ell}_x\times n_x}$ are model coefficients to be identified. Importantly, the number of terms in the function library $\b{\bar{\theta}_x}$ is at most, and in practice typically much smaller, than the number of terms in the function library $\b{\theta_x}$ used in manifold-aware identification, that is $\bar{\ell}_x\le\ell_x$. To achieve this, as a starting template, we build a function library $\b{\theta_z}(\b{x},\b{y})$ considering the full state $\b{z}$, including both slow and fast variables. Now, we assume that, based on physical priors for the system, the equations for the underlying dynamics have a known polynomial structure of degree $d_z$. Therefore, $\b{\theta_z}(\b{x},\b{y})$ is built including all monomials of up to that maximum degree $d_z$ in the state variables. This maximum monomial degree is typically rather low for the case of semi-discretised PDEs, including, for example, up to quadratic, cubic, and quartic terms for problems concerning incompressible fluid flows, structural dynamics with geometric nonlinearities, and heat transfer with radiation, respectively. Next, using a symbolic computation package, the learned model for the manifold in eq.~\eqref{ansatz_y} is substituted in place of all functions of the fast variables in the ansatz. Although this substitution results in a library that now only contains functions of $\b{x}$, as previously discussed, it will have linearly dependent terms, which is the root cause of ill-conditioning in the manifold-agnostic approach. However, using symbolic computation again, we may simplify and eliminate all redundant terms to obtain our manifold-informed function library used to identify the slow dynamics
\begin{equation}
    \b{\bar{\theta}_x{(\b{x})}_x} = \mathrm{Simplify}\left[\b{\theta_z}\left(\b{x},\b{\Xi_y}\T\b{\theta_y}(\b{x})\right) \right].
\end{equation}
If the manifold equation is sparse, then the result is a leaner and better conditioned function library that is consistent with our knowledge of the physics of the system and is informed by the structure of the manifold. In this work we use the SymPy Python library for symbolic computation, but other alternatives could be used as well. The approach used to build the manifold-informed library is explained in terms of pseudocode in algorithm~\ref{alg}.

\begin{algorithm}[H]
\caption{Manifold-informed library computation}\label{alg}
\begin{flushleft}
\textbf{Inputs:} \hspace{0.1 cm} data from slow variables $\b{X}$\\
\hspace{1.46 cm} SymPy symbols of the slow and fast variables, $\mathbf{x}$ and $\mathbf{y}$\\
\hspace{1.46 cm} degree of the manifold-agnostic library $d_z$\\
\hspace{1.46 cm} symbolic equation of the slow manifold $\mathbf{h}(\mathbf{x})$\\
\textbf{Outputs:} manifold-informed library $\b{\bar{\Theta}_x}$
\end{flushleft}
\begin{algorithmic}[1]
\State $\b{\theta_z}(\mathbf{x}, \mathbf{y}) \gets$ Generate symbolic library of degree $d_z$ using $(\mathbf{x}, \mathbf{y})$
\State $\b{\theta_z}(\mathbf{x}, \mathbf{h}(\mathbf{x})) \gets \b{\theta_z}(\mathbf{x}, \mathbf{y})$ \Comment{Substitute manifold equation via \texttt{SymPy}}
\State $d_y \gets$ highest polynomial degree in $\mathbf{h}(\mathbf{x})$
\State $d_x \gets d_z \cdot d_y$ \Comment{Maximum degree for the dynamics restricted to the slow
manifold}

\State Initialize empty list \texttt{updated\_functions} \Comment{Placeholder for $\b{\bar{\theta}_x}(\mathbf{x})$}
\For{each function $f_i$ in $\b{\theta_z}(\mathbf{x}, \mathbf{h}(\mathbf{x}))$}
    \State Decompose $f_i$ into monomial terms using \texttt{SymPy.expand()}
    \For{each term $t_j$ in $f_i$}
        \State Remove coefficient from $t_j$ \Comment{Retain only monomial structure}
        \State Add $t_j$ to \texttt{updated\_functions}
    \EndFor
\EndFor

\State \texttt{updated\_functions} $\gets$ Remove duplicates

\State $\b{\bar{\Theta}_x} \gets$ Evaluate each monomial in \texttt{updated\_functions} on data $\b{X}$

\State \textbf{return} $\b{\bar{\Theta}_x}$
\end{algorithmic}
\end{algorithm}

Now, equipped with our manifold-informed library, we assemble the data matrix
\begin{equation}
    \b{\bar{\Theta}_x}=\left[\b{\bar{\theta}_x}(\b{x}_1) \ \b{\bar{\theta}_x}(\b{x}_2) \ \cdots \ \b{\bar{\theta}_x}(\b{x}_m)\right]\T \in\mathbb{R}^{m\times \bar{\ell}_x}.\label{data_x}
\end{equation}
Lastly, the procedure concludes by using the STLS algorithm to approximate the solution to
\begin{equation}
  \b{\bar{\Xi}_x} = \underset{\b{\bar{\Xi}_x'}}{\mathrm{argmin}} \ \| \b{\dot{X}}  -\b{\bar{\Theta}_x} {\b{\bar{\Xi}_x'}} \|_F^2 + \bar{\lambda}_x \|\b{\bar{\Xi}_x'}\|_0,
\end{equation}
where the $\b{\dot{X}}$ data matrix is defined as in eq.~\eqref{data_x}, and $\bar{\lambda}_x$ is the sparsity promoting hyperparameter in this regression. As we show through our results in section~\S\ref{sec:results}, this method achieves robust, accurate, and fast identification of high-dimensional dynamics restricted to slow manifolds.

\subsection{Lower bound on the size of the manifold-informed library}

To quantify the potential reduction in library size with the proposed method, we present the following lower bound for the number of terms in the manifold-informed library

\begin{equation}
    \bar{\ell_x} \ge \sum_{j=0}^{d_z}\frac{(n_x+j-1)!}{j!(n_x-1)!} + d_z - \left\lfloor \frac{d_z}{d_y}\right\rfloor, \label{eq:bound}
\end{equation}
where $\lfloor \ \rfloor$ denotes the floor function that returns the integer that is closest and is less than or equal to its argument. The bound arises from counting the number of terms that result from the sparsest manifold equation resulting in the smallest achievable library. First, the summation accounts for all the terms in the library that are pure monomials in the slow variables $\b{x}$. Then we assume there is only one fast variable $y$. Moreover, in the best-possible scenario, a manifold equation of degree $d_y$ will depend on only one slow variable $x$ and will contain only one term, that is $y=x^{d_y}$. Substitution of this manifold equation into the function library for the slow variables will introduce at least $d_z - \left\lfloor d_z/d_y\right\rfloor$ new terms. These come from terms that were pure monomials in $y$ prior to the substitution, and subtracting those that result in monomials in $x$ of degree smaller or equal than $d_z$, which correspond to repeated terms. We point out that a tighter bound can probably be computed by accounting for the cross-terms between fast and slow variables, but we leave this for future efforts.

Comparing this lower bound with eq.~\eqref{eq:dense_library_size} for the manifold-aware library size, we see that the upper limit for the summation in the former is $d_x$, whereas in the latter it goes up to $d_z = d_x d_y$. Therefore, the reduction in library size can be significantly more impactful for higher $d_y$, that is, for more complex manifolds. Finally, we reiterate that, in addition to reducing computational cost and memory footprint, a leaner library also improves the conditioning of the subsequent regression problem.

\section{\label{sec:datasets}Numerical examples and dataset}

To demonstrate the effectiveness of our proposed manifold-informed identification, we have compiled a numerical dataset consisting of three slow-fast systems. One of these is the two-dimensional toy system introduced in section ~\S\ref{sec:methods2} and depicted in Fig.~\ref{fig:SINDyComparison}. The other two are high-dimensional systems arising from spatially discretized PDEs; these include a snap-through buckling beam and the fluid flow over a NACA 0012 airfoil. All the data is available on github.com/ben-herrmann.

\subsection{2D slow-fast system}

We consider the system introduced in section~\S\ref{sec:methods2} governed by eqs.~\eqref{eq:toy_system}. These equations were numerically integrated using SciPy's solve\_ivp over the time range $t\in[0, 3.8]$ with a sampling time step of $5\times10^{-3}$ yielding a total of $760$ data points per trajectory. Fifty trajectories were computed starting from random initial conditions $(x_0,y_0)$ close to the unstable equlibrium $(0,0)$ and on the slow manifold. The values of $x_0$ are sampled from a normal distribution with zero mean and variance $10^{-4}$. The corresponding initial values for $y_0$ were determined using eq.~\eqref{eq:toy_system DAE} for the manifold.

\subsection{Snap-through buckling beam}

\begin{figure}
    \centering
    \includegraphics[width=17.5cm]{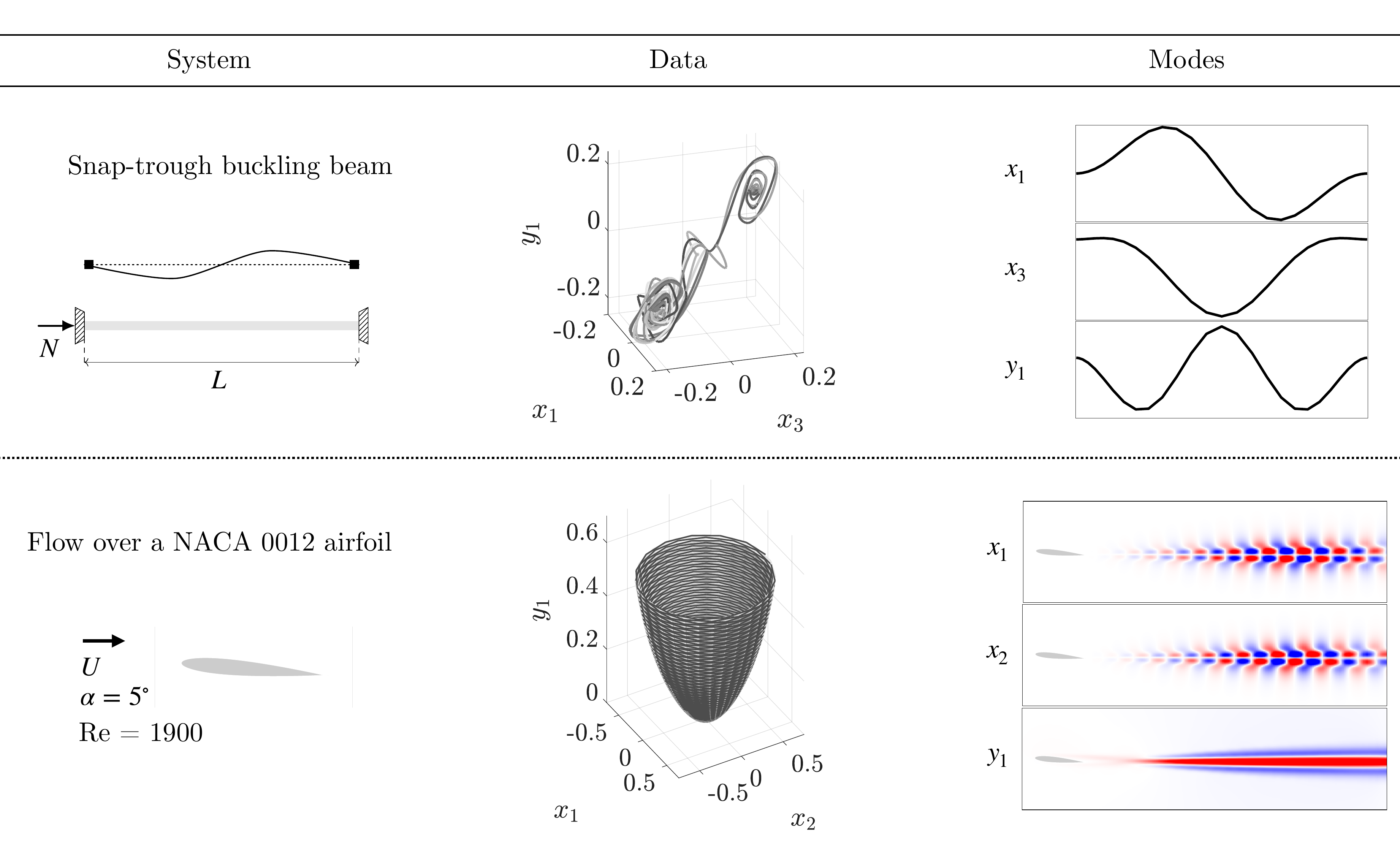}
    \caption{Schematics, the projection of trajectories to the latent space, and modes pertaining to both the snap-through buckling beam and flow through a NACA 0012 airfoil problems.}
    \label{fig:numerical examples}
\end{figure}

For our second example, we consider the dynamics of an axially loaded slender beam that is clamped at both of its edges and exhibits a snap-through buckling instability, the boundary conditions and load of this problem are displayed in Fig.~\ref{fig:numerical examples}. The model that we use for the beam dynamics is derived from Euler-Bernoulli beam theory incorporating a von Kármán geometric nonlinearity, and has been studied by several authors~\citep{wiebe2016prsa,zhong2021nd}. The dynamics of the system are governed by
\begin{equation}
    \label{eq:beam}
    \rho A\frac{\partial^2w}{\partial t^2} + C\frac{\partial w}{\partial t} + EI\frac{\partial ^4w}{\partial x^4} + \left(N - \frac{EA}{2L}\int_0^L\left(\frac{\partial w}{\partial x}\right) ^2dx\right)\frac{\partial ^2w}{\partial x^2}=0,
\end{equation}
where $w$ is the deflection of the beam, $L$ is its length, $t$ is time, $x$ is the spatial coordinate, $\rho$ and $ E$ represent the density and Young's modulus of the beam material, respectively, $A$ and $I$ are the cross-sectional area and the second moment of the area, $C$ accounts for linear viscous damping, and the term $N$ represents the axial loading. The dynamics are non-dimensionalized using $/2L$ and $L^2\sqrt{\rho A/EI}$ as length and time scales, respectively, and expressed as a system of first-order equations
\begin{subequations}
    \label{eq:beam_nd}
\begin{align}
    \frac{\partial w}{\partial t} &= v, \\
    \frac{\partial v}{\partial t} &= - c_1 v - \frac{\partial ^4w}{\partial x^4} - \left(c_2 - c_3\int_{-1}^1\left(\frac{\partial w}{\partial x}\right) ^2dx\right)\frac{\partial ^2w}{\partial x^2}=0,
\end{align}
\end{subequations}
where the dimensionless parameters are $c_1=CL^2/\sqrt{\rho A E I}$, $c_2=NL^2/EI$, and $c_3=AL^2/4I$. In this work, we fix these parameters at $c_1=2$, $c_2=30$ and $c_3=10^4$, for which the zero deflection configuration is linearly unstable and the system has two coexisting stable equilibria corresponding to the buckled upward and downward positions.

Equation~\eqref{eq:beam_nd} is spatially discretized using a Chebyshev pseudospectral collocation method~\citep{trefethen2000siam} with $31$ points and integrated in time using MATLAB. The clamped boundary conditions are imposed by making the deflection and its spatial derivative equal to zero at both ends of the beam. We assemble a dataset comprised of twenty trajectories generated using random small perturbations of the zero deflection configuration as initial conditions. The energy of these perturbations is normally distributed with zero mean and a standard deviation of $10^{-4}$. Each trajectory is integrated over $5$ time units, and the deflection and deflection velocity are registered every $0.02$ time units, yielding a total of $251$ snapshots.

The equations are also linearized about the $w=0$ equilibrium, and the discretized linear dynamics operator is used to perform a stability analysis. We compute and include in our dataset the eigenvalues and direct and adjoint eigenmodes obtained from the analysis. For the parameter values selected, the system has two unstable eigenvalues, corresponding to one even and one odd deflection mode. Moreover, the first two stable eigenmodes (in descending order according to their real parts) have the same shape as the unstable ones for the deflection, but the opposite sign in the deflection velocity. We take these four leading eigenmodes as our slow variables. Although the dynamics of the local unstable manifold of the equilibria can be characterized just using the unstable modes, to capture the stable (buckled) equilibria we need to include the first two stable eigenmodes as well. Therefore, the linear embedding matrices for the slow variables, $\b{V_x}$ and  $\b{W_x}$, are the first four direct and adjoint eigenmodes, respectively. For the fast variables, we find that including the subsequent $20$ eigenmodes in $\b{V_y}$ and  $\b{W_y}$ is enough to capture $99.3 \%$ of the sustained variance in the data. The deflection of a couple of slow modes and the first fast mode are shown in Fig.~\ref{fig:numerical examples}.

\subsection{Flow over a NACA 0012 airfoil}

For our last example, we consider a two-dimensional laminar flow over a NACA 0012 airfoil governed by the incompressible Navier-Stokes equations,

\begin{align}
    \rho\left(\frac{\partial\b{u}}{\partial t} + (\b{u}\cdot\nabla)\b{u}\right)&=-\nabla p+\mu\nabla^2\b{u},\label{eq:N-S}\\
     \nabla\cdot\b{u}&=0,
\end{align}

where $\b{u}$ represents the velocity field, $\rho$ and $\mu$ are the density and dynamic viscosity of the fluid, respectively, and $\nabla p$ stands for the pressure gradient. We choose a fixed angle of attack of $\alpha=5^{\circ}$ and a Reynolds number $\mathrm{Re}=1900$. The Reynolds number is a dimensionless quantity used in fluid mechanics that compares the relative importance of inertial and viscous effects and, for the flow over an airfoil, is defined as $\mathrm{Re}= \rho U c/\mu$, with $U$ and $c$ the free-stream velocity and the airfoil chord length, respectively. For these values of $\alpha$ and $\mathrm{Re}$, the system has a stable equilibrium corresponding to steady flow, but it is close to a supercritical Hopf bifurcation that occurs at $\mathrm{Re}_c\approx1950$ and leads to a limit cycle corresponding to periodic vortex shedding~\cite{gupta2023jfm}.

Direct numerical simulations (DNS) are performed using the spectral element code Nek5000~\cite{fischer2008nek5000}. We consider a rectangular computational domain extending from $-2c$ to $6c$ in the streamwise direction and $-2.5c$ to $2.5c$ in the crosswise direction, with the leading edge of the airfoil located at the origin. Spatial discretisation relies on a C-grid mesh embedded inside the rectangular domain, using $4368$ spectral elements with a polynomial order of $N = 5$. The same configuration was studied in~\citep{torres2024arxiv}.

In order to generate our initial conditions, we decide to compute the long-term solution for the flow at a higher Reynolds number, in this case after the bifurcation, to get a periodic vortex shedding solution for the system, and subsequently perform a DNS at a Reynolds number $\mathrm{Re}=1900$ with this solution as the initial condition. The periodic solution is computed by running DNS at ${Re}=2200$ from rest over a time horizon of $200$ time units. After that, $6$ DNS at a lower Reynolds number, $\mathrm{Re}=1900$, are computed for an additional $300$ time units. Given the periodic nature of the solution, we generate trajectories by selecting initial conditions that correspond to six different phases of the periodic vortex shedding. 

For the stability analysis we need to compute the equilibrium at $\mathrm{Re}=1900$ first. This is achieved with a long DNS because the solution at this Reynolds corresponds to a stable equilibrium. From here, we perform the stability analysis using an Arnoldi iteration, also implemented in nekStab~\cite{frantz2023amr}, with a Krylov basis dimension $m=350$ and a time step $\tau=0.1$. With this analysis, the leading $m=40$ eigenmodes and eigenvalues are computed, and then the two least stable modes are selected as the slow modes $\b{V_x}$. 

These modes are orthonormal, so we use $\b{W_x}=\b{V_x}$. To obtain an orthonormal basis for the fast modes $\b{V_y}=\b{W_y}$, we leverage the DNS data snapshots. Specifically, we perform POD on the velocity fluctuation data, aggregated from all trajectories, projected onto the orthogonal complement of the slow modes. That is, we do a singular value decomposition of $\left(\b{I}-\b{V_x}\b{V_x}\T\right)\b{X}$ and select the leading left singular vectors as the fast modes, employing a cutoff of $99\%$ of the singular value energy. As a result, we end up with two slow modes and seven fast modes, all with a dimension of $n_q=314,496$. The streamwise component of the velocity field for the slow modes, as well as that of the first fast mode, is depicted in Fig.~\ref{fig:numerical examples}.

\section{\label{sec:results}Results and discussion}

In this section, we apply our manifold-informed identification on the numerical examples described in the previous section, comparing its performance to both manifold-agnostic and manifold-aware approaches. For each example, half of the trajectories are randomly selected as training data, while the other half is reserved for validation. We use the PySINDy Python package to perform the SINDy regressions~\cite{desilva2020joss, kaptanoglu2022joss}. In particular, we employ second-order central differences as the differentiation scheme and the STLS algorithm with normalization as the optimizer; where each column of a matrix library is divided by its norm. Regarding the hyperparameter $\lambda$ in the regressions, we performed a sweep across a logarithmically spaced range of $100$ points, spanning from $10^{-3}$ to $1$, for each optimization problem.

Once the governing equations have been identified, we use the initial conditions from the validation data to integrate the equations forward in time. This process allows us to compute the reconstructed slow variables, represented as $\b{\hat{x}}$ and the fast variables, denoted as $\b{\hat{y}}$. If necessary, the high-dimensional space of the system can also be reconstructed using the equation
\begin{equation*}
    \b{\hat{q}}=\b{V_x \hat{x}} +\b{V_y \hat{y}}.
\end{equation*}

For manifold-agnostic libraries, the reconstruction procedure is straightforward. After identifying the differential equations for both sets of variables, they can be integrated to obtain $\b{\hat{x}}$ and $\b{\hat{y}}$. In contrast, utilising manifold-informed (manifold-aware) libraries involves additional steps. Initially, we derive the equation of the slow manifold. Subsequently, the manifold-informed (manifold-aware) library is assembled, and the differential equations for the slow variables are determined. These equations are then integrated to produce $\b{\hat{x}}$. Finally, the computed trajectory $\b{\hat{x}}$ is input into the manifold equation to calculate $\b{\hat{y}}$.

A metric to evaluate concerning library matrices is their condition number, denoted as $\kappa$. In the context of numerical regressions, the condition number serves as an indicator of the sensitivity of the solution to noisy data. It significantly influences the accuracy of the computed pseudoinverse matrices, which are essential for solving these optimisation problems. Particularly, in a scenario involving a SINDy problem that encompasses two different library matrices, the matrix with a higher condition number is deemed to be the worse poised for the algorithm.

\subsection{Results for the 2D slow-fast system }

For this dynamical system, as presented in Fig.~\ref{fig:SINDyComparison}, the degree of the system equations is $d_z=2$. As previously stated, the slow manifold of the system can be approximated using a polynomial series of order $d_y=4$; then, the maximum degree for both the aware and informed libraries is $d_x=8$. The identified equations by the manifold-agnostic approach are given by
\begin{subequations}
\begin{align}
    \dot{x}&=1.000x-1.000xy,\\
    \dot{y}&=0.649x^2,
\end{align}
\end{subequations}
the condition number of its library matrix is, $\kappa(\b{\Theta_z})=2.54\times10^4$. Whereas the identified manifold equation of both manifold-aware and manifold-informed methods is expressed by
\begin{equation}
    y = 0.331x^2 + 0.056x^4,
\end{equation}
meanwhile, the slow variable differential equation derived from a manifold-aware method with $\kappa(\b{\Theta_x})=9.85\times10^2$, which is the same as the one derived from a manifold-informed procedure with $\kappa(\b{\bar{\Theta}_x})=6.58\times10^2$ is
\begin{equation}
    \dot{x} = 1.000 x -0.328x^3 -0.063x^5.
\end{equation}
Due to the simplicity of the equations involved in this system, we can calculate the model identification error with
\begin{equation}
    e_{\b{\Xi}}=\frac{\|\b{\Xi}-\b{\Xi}_{\text{true}}\|_F}{\|\b{\Xi}_{\text{true}}\|_F},
\end{equation}
for the manifold-agnostic identification, the model error is 56.57\%. For the manifold equation discovery, which is common for the other methods, the error model is 3.51\%, the error for the aware technique is 1.83\%, while our informed procedure shows an error of 1.83\%. Yet, the reconstruction error is $0.23\%$ for the agnostic method, $0.06\%$ for the aware one, and $0.03\%$ for the informed approach.

\subsection{Results for high-dimensional systems}

As previously explained, the initial phase in the assembly of a manifold-informed library involves the sparse identification of the underlying slow manifold. In the buckling beam example, a polynomial library of degree $d_y=2$ yields a manifold equation that replicates the validation data with an error of $4.3\%$. Conversely, for the vortex shedding case, the error is of $1.4\%$, but this is achieved with a manifold of degree $d_y=4$. The reconstruction errors of the fast variables are quantified via the normalized mean trajectory error (NMTE) defined in ~\citep{cenedese2022natcomms, axaas2023nd_a} as
\begin{equation}
\label{eq:NTME}
\text{NTME} = \frac{1}{m}\frac{1}{\underset{j}{\max}\|\b{y}(t_j)\|_{\b{M}}}\sum_{j=1}^m\|\b{y}(t_j)-\b{\hat{y}}(t_j)\|_{\b{M}},
\end{equation}
where the energy norm for $\b{x}\in\mathbb{R}^n$ is given by $\|\b{x}\|_{\b{M}}=(\b{x}^T \b{M}\b{x})^{\frac{1}{2}}$, with $\b{M}\in\mathbb{R}^{n\times n}$ being a symmetric positive definite matrix. For these datasets in particular, $\b{M}$ is the mass matrix of the dynamic systems in question.

Subsequently, leveraging the known polynomial structure of the physical examples, we can ascertain the maximum degree of the dynamics on the slow manifold. Specifically, in the buckling beam problem, the physics degree corresponds to $d_z=3$, which means that $d_x=6$. In contrast, for the fluid flow example $d_z=2$, resulting in $d_x=8$.

The resultant sizes of the manifold-informed libraries, denoted as $\bar{\ell}_x$, are shown in Table~\ref{tab:conditions and sizes} along with the library sizes for both manifold-agnostic and manifold-aware approaches. The quantity $\bar{\ell}_x$ is $135$ for the buckling beam problem and $40$ for the fluid flow problem. For the buckling beam example, this reflects a substantial reduction of $95.4\%$ relative to the manifold-agnostic library size $\ell_z$, and a $35.4\%$ decrease when compared to the manifold-aware library size $\ell_x$. Concerning the airfoil problem, the reductions in library size are quantified at $25.9\%$ for $\ell_z$ and $9.1\%$ for $\ell_x$. The manifold-informed library for the buckling beam effectively reduces the quantity of monomials for each degree exceeding $d_z$ in comparison to the manifold-aware library, whereas for the fluid flow case, it only eliminates terms of degree eight.

\begin{table}[h]
\begin{tabular}{c|c|c|c|c|c|c}
\hline
System                     & $\ell_z$   & $\kappa(\b{\Theta_z})$       & $\ell_x$ & $\kappa(\b{\Theta_x})$         & $\bar{\ell}_x$  & $\kappa(\b{\bar{\Theta}_x})$       \\ \hline
Snap-through buckling beam & $2924$ & $1.44\times10^{17}$ & $209$ & $1.23\times10^{9}$ & $135$ & $3.30\times10^{7}$ \\
NACA 0012 airfoil          & $54$   & $1.39\times10^5$  & $44$  & $7.80\times10^3$  & $40$  & $7.17\times10^3$  \\ \hline
\end{tabular}
\caption{Comparison of the condition number and library sizes between the manifold-agnostic, manifold-aware and manifold-informed library matrices on numerical experiments.}
\label{tab:conditions and sizes}
\end{table}

In the buckling beam example, as displayed in Table~\ref{tab:conditions and sizes}, the condition number of the manifold-informed library exhibits a reduction of ten orders of magnitude compared to that of the agnostic library, and is two orders of magnitude lower than the condition number of the aware library. Meanwhile, for the vortex shedding problem, the condition number differs by two orders of magnitude in relation to the manifold-agnostic library and by less than one order of magnitude compared to the manifold-aware library.

\begin{figure}
    \centering
    \includegraphics[width=17.5cm]{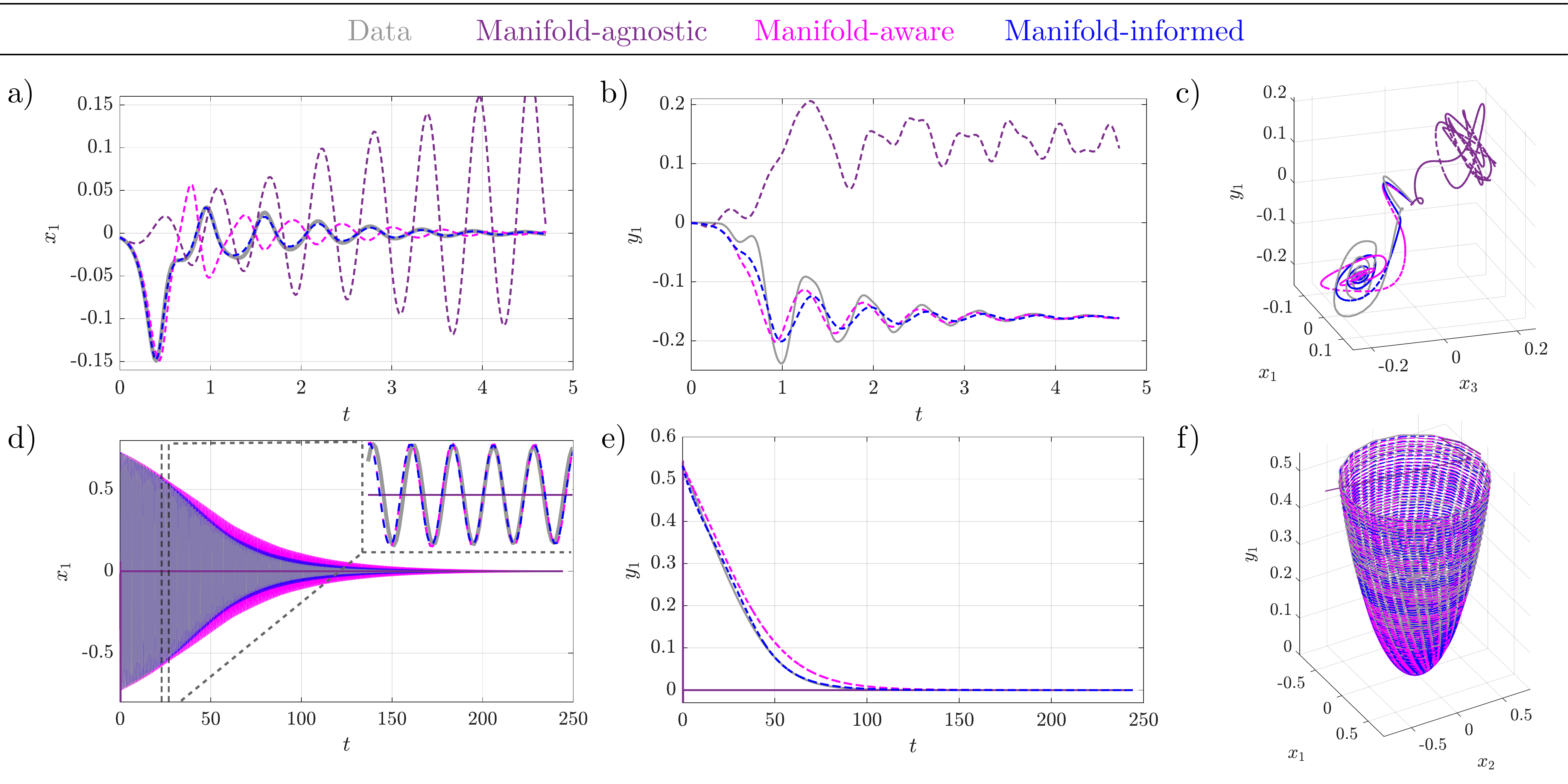}
    \caption{Reconstruction of validation trajectories using the equations identified with manifold-agnostic, manifold-aware and manifold-informed methods. Subfigures (a) and (b) show the slow and fast variables of the buckling beam problem. (c) Depiction of the phase space of the buckling beam problem. Subfigures (d) and (e) show the slow and fast variables of the vortex shedding problem. (f) Representation of the phase space of the vortex shedding problem.}
    \label{fig:Coeffs}
\end{figure}

After integrating the selected models for each library type, we can compare them to assess how effectively each model reconstructs the validation data. Notably, the time series of the slow and fast variables as depicted in Fig.~\ref{fig:Coeffs} reveal that the outputs derived from our manifold-informed library exhibit a temporal evolution more closely aligned with the validation data compared to those from both the manifold-aware and agnostic libraries.

\begin{figure}
    \centering
    \includegraphics[width=17.5cm]{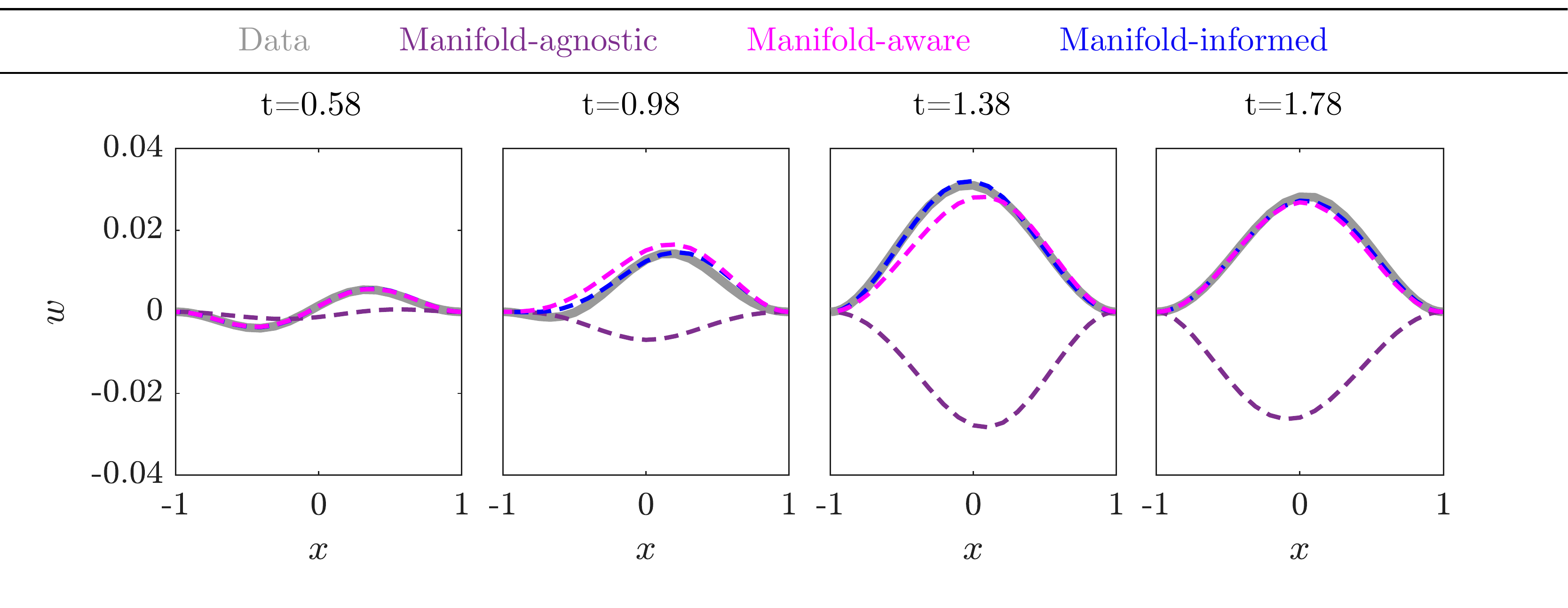}
    \caption{Snapshots depicting the reconstruction of the deflection in the physical space of one validation trajectory of the snap-trough buckling beam.}
    \label{fig:BeamSnapshots}
\end{figure}

\begin{figure}
    \centering
    \includegraphics[width=17.5cm]{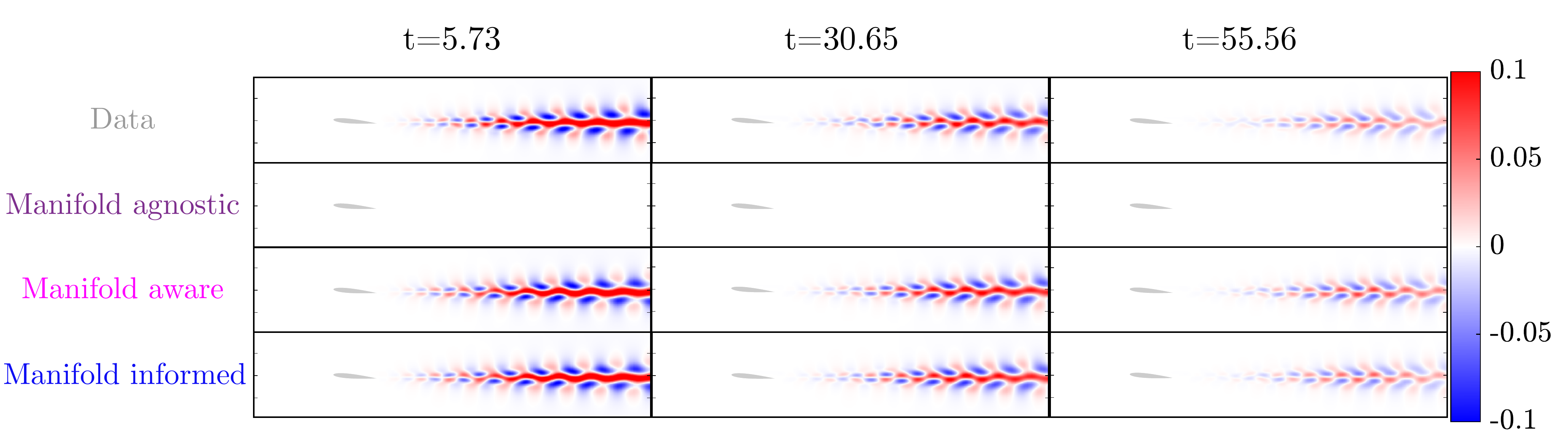}
    \caption{Snapshots depicting the reconstruction of the streamwise velocity in the physical space of one validation trajectory of the vortex shedding on a flow through an airfoil.}
    \label{fig:AirfoilSnapshots}
\end{figure}

Consistent with the observations regarding slow and fast variables, the reconstruction of the validation data in physical space is also more accurate for results obtained with the manifold-informed model than for those obtained with manifold-aware and agnostic approaches. This is particularly evident in the qualitative evaluations represented in snapshots from  Fig.~\ref{fig:BeamSnapshots} and Fig.~\ref{fig:AirfoilSnapshots}. Specifically, in the buckling beam problem, the manifold-informed library consistently produces reconstructions that closely match the validation data, whereas the manifold-aware library shows slight deviations. Conversely, reconstructions produced from the manifold-agnostic library initially align closely with the data but eventually deviate more significantly over time than those derived from the aware approach; furthermore, for certain trajectories, the agnostic reconstructions converge into an equilibrium that differs from that of the original data. For the vortex shedding case, both the manifold-informed and manifold-aware models effectively capture the decay rate of oscillations in the flow, although the amplitude of the oscillations reconstructed using our method is more closely aligned with the amplitude of the data, while the agnostic case predominantly reflects the base flow regime.

To quantitatively assess the fidelity of these reconstructions in physical space, we employ the NMTE defined in eq.~\eqref{eq:NTME} but for $\b{q}$ instead of $\b{y}$. The error metrics documented in Table~\ref{tab: NTME} consistently indicate that our proposed method outperforms both the agnostic and aware approaches. Errors exceeding $>100\%$ are shown for cases where unstable models are identified and the solution blows up. These results align with the hypothesis regarding the ill-conditioning of the library matrices.

\begin{table}[h!]
\begin{tabular}{c|ccc}
\hline
\multirow{2}{*}{System}    & \multicolumn{3}{c}{NTME}                                                                             \\ \cline{2-4} 
                           & \multicolumn{1}{c|}{Manifold-agnostic}     & \multicolumn{1}{c|}{Manifold-aware} & Manifold-informed \\ \hline
Snap-through buckling beam & \multicolumn{1}{c|}{$96.5\pm41.8$ \%}              & \multicolumn{1}{c|}{$10.9\pm4.3$ \%}           & $6.8\pm0.7$ \%               \\
NACA 0012 airfoil          & \multicolumn{1}{c|}{$>100$ \%} & \multicolumn{1}{c|}{$0.6\pm0.1$ \%}        & $0.5\pm0.1$ \%               \\ \hline
\end{tabular}
\caption{Normalized mean trajectory error of the reconstructed trajectories $\b{\hat{q}}$ derived from a manifold-informed, manifold-aware, and manifold-agnostic libraries.}
\label{tab: NTME}
\end{table}

\section{\label{sec:conclusions}Conclusions}

We developed an extension of the sparse identification of nonlinear dynamics (SINDy) algorithm specifically tailored for systems exhibiting slow-fast dynamics. By identifying, and subsequently leveraging, the structure of the underlying slow manifold, our method constructs a manifold-informed function library that avoids the redundancy and ill-conditioning inherent in manifold-agnostic identification approaches. At the same time, our approach reduces the size of the function library compared to the library that would be used in existing manifold-aware methods, such as SINDy with autoencoders, operator inference on manifolds, and data-driven SSM identification. This is achieved by only including nonlinearities that are consistent with both, the previously known structure of the physical system and the manifold equation learned from data. Moreover, we provide a lower bound for the achievable reduction in library size that allows assessing the potential impact of the approach for a particular scenario.

We demonstrate our method on numerical examples, including a snap-through buckling beam and flow over a NACA0012 airfoil, showcasing its ability to extract interpretable and accurate slow manifold dynamics. Compared to previous work, our approach extends the applicability of SINDy to slow-fast systems while addressing key limitations in function library size and regression ill-conditioning. We find that our method consistently achieves lower model identification errors and improved reconstruction accuracy compared to both manifold-agnostic and manifold-aware approaches. Notably, our approach leads to leaner and better conditioned function library matrices, thus enhancing numerical stability and computational efficiency.

The proposed method can be applied to data from any slow manifold, such as stable, center, and unstable manifolds, as well as SSMs attached to steady states. In the latter case, our approach can be easily integrated into the recent fast data-driven SSM identification~\citep{axaas2023nd_a} to reduce the size of the polynomial library that is used for the dynamics restricted to the SSM, which may be the bottleneck for highly nonlinear manifolds. Therefore, using a manifold-informed library can potentially enable scaling of fast data-driven SSM identification to dynamics on higher-dimensional manifolds.

The developed method can be easily coupled with other available SINDy extensions. Future work will focus on extending this methodology to parameterized systems, where both the dynamics and the manifold may change with parameters.

\begin{acknowledgments}
We gratefully acknowledge J. Lemus and N. Torres for their helpful comments and insightful discussions. We also remark that part of this work was developed while B. Herrmann was affiliated to the Department of Mechanical Engineering at Universidad de Chile. This work was funded by ANID Fondecyt 1250693.
\end{acknowledgments}

\appendix

\section{Appendixes}\label{app1}

\subsection{Derivation of the slow manifold equation}
Consider a 2D dynamical system with an equilibrium point at $\b{x}=\b{0}$. We suppose that the manifold in the neighborhood of the equilibrium can be approximated by a series of polynomials
\begin{equation*}
    h(x)=ax^2+bx^3+cx^4+\mathcal{O}(5)\quad\quad a,b,c\in\mathbb{R},
\end{equation*}
this manifold must be in consistent with the tangency condition
\begin{equation*}
    \dot{y}=\frac{\partial h}{\partial x}\dot{x}.
\end{equation*}
Consider the slow-fast system defined by the following equations
\begin{align*}
    \dot{x}&=x-xy,\\
    \dot{y}&=-y+x^2,
\end{align*}
now, we apply the tangency condition to determine the coefficients of the manifold approximation
\begin{align*}
    \dot{y}&=-(ax^2+bx^3+cx^4)+x^2,\\
    \frac{\partial h}{\partial x}&=2ax + 3bx^2 + 4cx^3 +\mathcal{O}(4),\\
    \dot{x}&=x-x(ax^2+bx^3+cx^4+\mathcal{O}(4)),\\
    &\Longrightarrow 2ax^2 + 3bx^3 + (4c-2a^2)x^4-(1-a)x^2+bx^3+cx^4=0,\\
    &\Longrightarrow a=1/3\quad\wedge\quad b=0\quad\wedge\quad c=\frac{2}{45},\\
    &\Longrightarrow h(x)\simeq\frac{x^2}{3}+\frac{2x^4}{45}.
\end{align*}

\subsection{The relevance of the fast modes}

\begin{figure}[H]
    \centering
    \includegraphics[width=17.5cm]{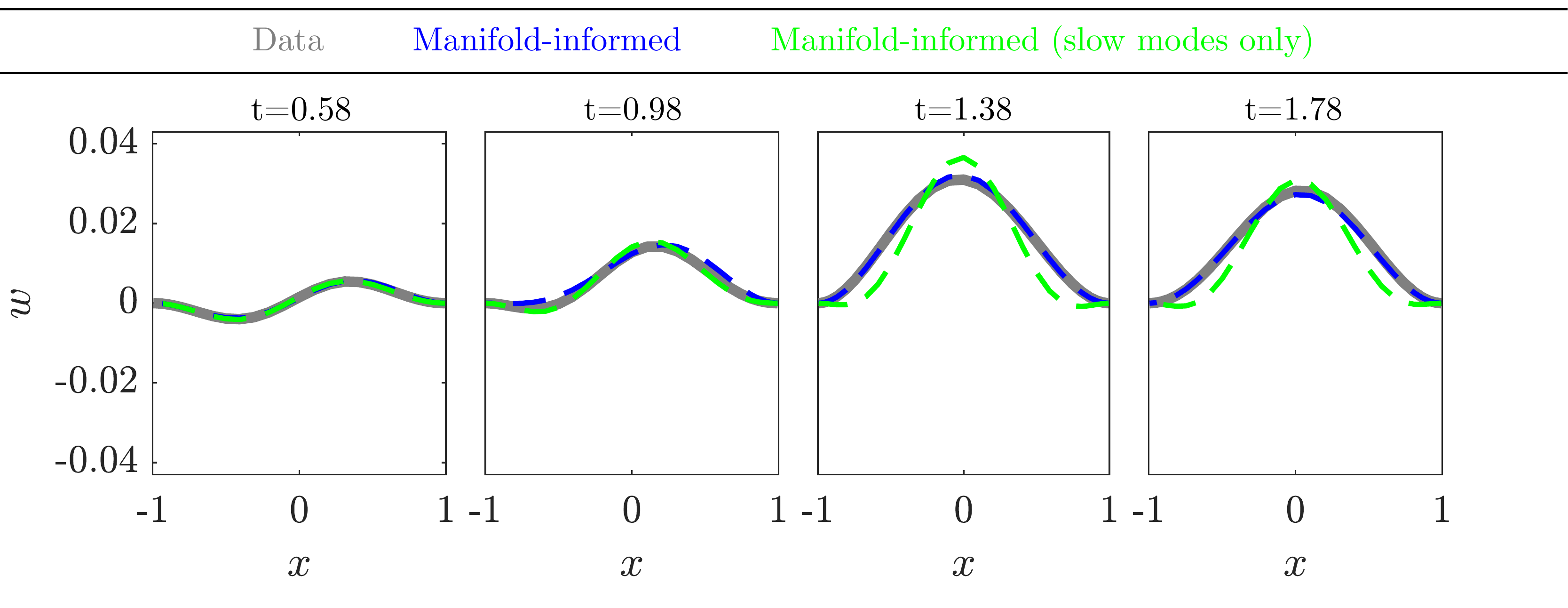}
    \caption{Snapshots of the reconstruction in the physical space for one of the validation trajectories of the snap-through buckling beam. It compares the results of using a linear embedding that only contains slow modes and a linear embedding that incorporates both slow and fast modes.}
    \label{fig:BeamSnapshotsSlow}
\end{figure}

Because the dynamics of a slow-fast system can be written into a DAE depending only on the slow variables, one may raise the question of whether the fast variables are imperative for a good enough reconstruction of a high-dimensional physical space. As illustrated in Fig.~\ref{fig:BeamSnapshotsSlow}, when $\b{W_y}=\b{V_y}=\b{0}$ (a simplification for clarity), the reconstruction that excludes the fast modes is noticeably inferior to the reconstruction that includes all the modes.

\providecommand{\noopsort}[1]{}\providecommand{\singleletter}[1]{#1}%

\end{document}